\documentclass{amsart}
\usepackage{amssymb}
\usepackage{amsmath}
\usepackage{amsfonts}

\theoremstyle{plain}

\newtheorem{lemma}{Lemma}

\numberwithin{equation}{section}
\input{tcilatex}

\begin{document}
\title[Meridian Surface of Weingarten Type]{Meridian Surface of Weingarten
Type in 4-dimensional Euclidean Space $\mathbb{E}^{4}$}
\author{G\"{u}nay \"{O}zt\"{u}rk}
\address{Kocaeli University, Department of Mathematics\\
41380 Kocaeli, TURKEY}
\email{ogunay@kocaeli.edu.tr}
\author{Bet\"{u}l Bulca}
\address{Uluda\u{g} University, Department of Mathematics\\
16059 Bursa, TURKEY}
\email{bbulca@uludag.edu.tr}
\author{Beng\"{u} K\i l\i \c{c} Bayram}
\address{Bal\i kesirUniversity, Department of Mathematics\\
Bal\i kesir, TURKEY}
\email{genguk@bal\i kesir.edu.tr}
\author{Kadri Arslan}
\address{Uluda\u{g} University, Department of Mathematics\\
16059 Bursa, TURKEY}
\email{arslan@uludag.edu.tr}

\begin{abstract}
In this paper, we study meridian surfaces of Weingarten type in Euclidean
4-space $\mathbb{E}^{4}.$ We give the neccessary and sufficient conditions
for a meridian surface in $\mathbb{E}^{4}$ to become \ Weingarten type.
\end{abstract}

\maketitle

\section{Introduction}

A surface $M$ in $\mathbb{E}^{n}$ is called Weingarten surface if there
exist a non-trivial function%
\begin{equation}
\Psi (K,H)=0  \label{A1}
\end{equation}%
between the Gauss curvature $K$ and mean curvature $H$ of the surface $M$.
The existence of a non-trivial functional relation $\Psi (K,H)=0$ on a
surface $M$ parametrized by a patch $X(u,v)$ is equivalent to the vanishing
of the corresponding Jacobian determinant, namely%
\begin{equation}
\left \vert \frac{\partial (K,H)}{\partial (u,v)}\right \vert =0.  \label{A2}
\end{equation}%
The condition (\ref{A2}) that must be satisfied for the Weingarten surface $%
M $ leads to 
\begin{equation}
K_{u}H_{v}-K_{v}H_{u}=0  \label{A3}
\end{equation}%
with subscripts denoting partial derivatives.

For the study of these surfaces, W. K\"{u}hnel \cite{K} investigated ruled
Weingarten surface in a Euclidean 3-space $\mathbb{E}^{3}$. Further, D. W.
Yoon \cite{Y} classified ruled linear Weingarten surface in $\mathbb{E}^{3}.$
Meanwhile, F. Dillen and W. K\"{u}hnel \cite{DK} and Y. H. Kim and D. W.
Yoon \cite{KY} gave a classification of ruled Weingarten surfaces in a
Minkowski 3-space $\mathbb{E}_{1}^{3}$. . Recently, M. I. Munteanu and I.
Nistor \cite{MN} and R. Lopez (\cite{L1,L2}) studied polynomial translation
Weingarten surfaces in a Euclidean 3-space.

The study of meridian surfaces in $\mathbb{E}^{4}$ was first introduced by
G. Ganchev and V. Milousheva (See, \cite{GM} and \cite{GM2}). They construct
a surface $M^{2}$ in $\mathbb{E}^{4}$ in the following way: 
\begin{equation}
M^{2}:X(u,v)=f(u)\,r(v)+g(u)\,e_{4},\quad u\in I,\,v\in J  \label{A4}
\end{equation}%
where $f=f(u),\, \,g=g(u)$ are non-zero smooth functions, defined in an
interval $I\subset \mathbb{R}$, such that $(f^{\prime }(u))^{2}+(g^{\prime
}(u))^{2}=1,\, \,u\in I$ and $r=r(v)$ $\, \left( v\in J\subset \mathbb{R}%
\right) $ is a curve on $S^{2}(1)$ parameterized by the arc-length and $%
e_{4} $ is the fourth vector of the standard orthonormal frame in $\mathbb{E}%
^{4}.$ See also \cite{BAM} for the classification of meridian surfaces in
4-dimensional Euclidean space $\mathbb{E}^{4}$ which have pointwise 1-type
Gauss map.

In this paper, we study meridian surfaces of Weingarten type in
4-dimensional Euclidean space $\mathbb{E}^{4}.$ We proved the following
theorem:

\textbf{Main Theorem.} Let $M^{2}$ be a meridian surface given with the
parametrization (\ref{C2}). Then $M^{2}$ is a Weingarten surface if and only
if \ $M^{2}$ is one of the following surfaces;

$i\mathbf{)}$ a planar surface lying in the constant $3$-dimensional space
spanned by $\{x,y,n_{2}\}$,

$ii)$ a developable ruled surface in a $3$-dimensional Euclidean space $%
\mathbb{E}^{3}$,

$iii)$ a developable ruled surface in a $4$-dimensional Euclidean space $%
\mathbb{E}^{4},$

$iv\mathbf{)}$ a surface given with the surface patch%
\begin{eqnarray*}
X(u,v) &=&\left( \frac{\cos \left( au+ac_{1}\right) }{a}+c_{2}\right) \,r(v)+
\\
&&+\left( \frac{2\left( \sin \left( au+ac_{1}\right) -1\right) \sqrt{1+\sin
\left( au+ac_{1}\right) }}{\cos \left( au+ac_{1}\right) }\right) e_{4},
\end{eqnarray*}

$iv\mathbf{)}$ a surface given with the surface patch%
\begin{eqnarray*}
X(u,v) &=&\pm \frac{1}{2}\left( \sqrt{-\frac{a}{\left( e^{\frac{u}{b}%
}\right) ^{2}\left( e^{\frac{c}{b}}\right) ^{2}}}\left( \left( e^{\frac{u}{b}%
}\right) ^{2}\left( e^{\frac{c}{b}}\right) ^{2}+1\right) \right) \,r(v) \\
&&\pm \frac{1}{2}\left( \sqrt{\frac{4b^{2}\left( e^{\frac{u}{b}}\right)
^{2}\left( e^{\frac{c}{b}}\right) ^{2}+a\left( e^{\frac{u}{b}}\right)
^{4}\left( e^{\frac{c}{b}}\right) ^{4}-2a\left( e^{\frac{u}{b}}\right)
^{2}\left( e^{\frac{c}{b}}\right) ^{2}+a}{b^{2}\left( e^{\frac{u}{b}}\right)
^{2}\left( e^{\frac{c}{b}}\right) ^{2}}}\right) e_{4}
\end{eqnarray*}

where $a,$ $b,$ $c,$ $c_{1},$ $c_{2}$ are real constants.

\section{Basic Concepts}

Let $M$ be a smooth surface in $\mathbb{E}^{n}$ given with the patch $X(u,v)$
: $(u,v)\in D\subset \mathbb{E}^{2}$. The tangent space to $M$ at an
arbitrary point $p=X(u,v)$ of $M$ span $\left \{ X_{u},X_{v}\right \} $. In
the chart $(u,v)$ the coefficients of the first fundamental form of $M$ are
given by 
\begin{equation}
E=\left \langle X_{u},X_{u}\right \rangle ,F=\left \langle X_{u},X_{v}\right
\rangle ,G=\left \langle X_{v},X_{v}\right \rangle ,  \label{B1}
\end{equation}%
where $\left \langle ,\right \rangle $ is the Euclidean inner product. We
assume that $W^{2}=EG-F^{2}\neq 0,$ i.e. the surface patch $X(u,v)$ is
regular.\ For each $p\in M$, consider the decomposition $T_{p}\mathbb{E}%
^{n}=T_{p}M\oplus T_{p}^{\perp }M$ where $T_{p}^{\perp }M$ is the orthogonal
component of $T_{p}M$ in $\mathbb{E}^{n}$.

Let $\chi (M)$ and $\chi ^{\perp }(M)$ be the space of the smooth vector
fields tangent to $M$ and the space of the smooth vector fields normal to $M$%
, respectively. Given any local vector fields $X_{1},$ $X_{2}$ tangent to $M$%
, consider the second fundamental map $h:\chi (M)\times \chi (M)\rightarrow
\chi ^{\perp }(M);$%
\begin{equation}
h(X_{i},X_{_{j}})=\widetilde{\nabla }_{X_{_{i}}}X_{_{j}}-\nabla
_{X_{_{i}}}X_{_{j}}\text{ \  \  \ }1\leq i,j\leq 2.  \label{B2}
\end{equation}%
where $\nabla $ and $\overset{\sim }{\nabla }$ are the induced connection of 
$M$ and the Riemannian connection of $\mathbb{E}^{n}$, respectively. This
map is well-defined, symmetric and bilinear.

For any arbitrary orthonormal frame field $\left \{
N_{1},N_{2},...,N_{n-2}\right \} $ of $M$, recall the shape operator $A:\chi
^{\perp }(M)\times \chi (M)\rightarrow \chi (M);$%
\begin{equation}
A_{N_{k}}X_{j}=-(\widetilde{\nabla }_{X_{j}}N_{k})^{T},\text{ \  \  \ }%
X_{j}\in \chi (M).  \label{B3}
\end{equation}%
This operator is bilinear, self-adjoint and satisfies the following equation:%
\begin{equation}
\left \langle A_{N_{k}}X_{j},X_{i}\right \rangle =\left \langle
h(X_{i},X_{j}),N_{k}\right \rangle =c_{ij}^{k}\text{, }1\leq i,j\leq 2;\text{
}1\leq k\leq n-2  \label{B4}
\end{equation}%
where $c_{ij}^{k}$ are the coefficients of the second fundamental form.

The equation (\ref{B2}) is called Gaussian formula, and%
\begin{equation}
h(X_{i},X_{j})=\overset{n-2}{\underset{k=1}{\sum }}c_{ij}^{k}N_{k},\  \  \  \  \
1\leq i,j\leq 2.\text{ }  \label{B5}
\end{equation}

Then the Gauss curvature $K$ of a regular patch $X(u,v)$ is given by

\begin{equation}
K=\frac{1}{W^{2}}\sum%
\limits_{k=1}^{n-2}(c_{11}^{k}c_{22}^{k}-(c_{12}^{k})^{2}).  \label{B6}
\end{equation}

Further, the mean curvature vector of a regular patch $X(u,v)$ is given by 
\begin{equation}
\overrightarrow{H}=\frac{1}{2W^{2}}%
\sum_{k=1}^{n-2}(c_{11}^{k}G+c_{22}^{k}E-2c_{12}^{k}F)N_{k}.  \label{B7}
\end{equation}%
where $E,F,G$ are the coefficients of the first fundamental form and $%
c_{ij}^{k}$ are the coefficients of the second fundamental form.

The norm of the mean curvature vector $H=\left \Vert \overrightarrow{H}%
\right \Vert $ is called the mean curvature of $M$. The mean curvature $H$
and the Gauss curvature $K$ play the most important roles in differential
geometry for surfaces \cite{Ch2}. Recall that a surface $M$ is said to be 
\textit{flat} (resp. \textit{minimal})\textit{\ } if its Gauss curvature
(resp. mean curvature vector) vanishes identically \cite{Ch}.

\section{Meridian Surfaces in $\mathbb{E}^{4}$}

Let $\{e_{1},e_{2},e_{3},e_{4}\}$ be the standard orthonormal frame in $%
\mathbb{E}^{4}$, and $S^{2}(1)$ be a 2-dimensional sphere in $\mathbb{E}%
^{3}=span\{e_{1},e_{2},e_{3}\}$, centered at the origin $O$. We consider a
smooth curve $c:r=r(v),\,v\in J,\, \,J\subset \mathbb{R}$ on $S^{2}(1)$,
parameterized by the arc-length ($r^{\prime }{}^{2}(v)=1$). We denote $%
t(v)=r^{\prime }(v)$ and consider the moving frame field $\{t(v),n(v),r(v)\}$
of the curve $c$ on $S^{2}(1)$. With respect to this orthonormal frame field
the following Frenet formulas hold good: 
\begin{equation}
\begin{array}{l}
\vspace{2mm}r^{\prime }(v)=t(v); \\ 
\vspace{2mm}t^{\prime }(v)=\kappa (v)\,n(v)-r(v); \\ 
\vspace{2mm}n^{\prime }(v)=-\kappa \,(v)t(v),%
\end{array}
\label{C1}
\end{equation}%
where $\kappa $ is the spherical curvature of $c$.

Let $f=f(u),\, \,g=g(u)$ be non-zero smooth functions, defined in an
interval $I\subset \mathbb{R}$, such that $(f^{\prime }(u))^{2}+(g^{\prime
}(u))^{2}=1,\, \,u\in I$. Now we construct a surface $M^{2}$ in $\mathbb{E}%
^{4}$ in the following way: 
\begin{equation}
M^{2}:X(u,v)=f(u)\,r(v)+g(u)\,e_{4},\quad u\in I,\,v\in J  \label{C2}
\end{equation}

The surface $M^{2}$ lies on the rotational hypersurface $M^{3}$ in $\mathbb{E%
}^{4}$ obtained by the rotation of the meridian curve $\alpha :u\rightarrow
(f(u),g(u))$ around the $Oe_{4}$-axis in $\mathbb{E}^{4}$. Since $M^{2}$
consists of meridians of $M^{3}$, we call $M^{2}$ a \textit{meridian surface 
}(see, \cite{GM}).

The tangent space of $M^{2}$ is spanned by the vector fields: 
\begin{equation}
\begin{array}{l}
\vspace{2mm}X_{u}(u,v)=f^{\prime }(u)r(v)+g^{\prime }(u)e_{4}; \\ 
\vspace{2mm}X_{v}(u,v)=f(u)\,t(v),%
\end{array}
\label{C3}
\end{equation}%
and hence the coefficients of the first fundamental form of $M^{2}$ are $%
E=1;\, \,F=0;\, \,G=f^{2}(u)$. Without louse of generality we can take $%
g^{\prime }(u)\neq 0$. Taking into account (\ref{C1}), we calculate the
second partial derivatives of $X(u,v)$: \ 
\begin{equation}
\begin{array}{l}
\vspace{2mm}X_{uu}(u,v)=f^{\prime \prime }(u)r(v)+g^{\prime \prime
}\,(u)e_{4}; \\ 
\vspace{2mm}X_{uv}(u,v)=f^{\prime }(u)t(v); \\ 
\vspace{2mm}X_{vv}(u,v)=f(u)\kappa (v)\,n(v)-f(u)\,r(v).%
\end{array}
\label{C31}
\end{equation}

Let us denote $X=X_{u},\, \,Y={\frac{X_{v}}{f}=t}$ and consider the
following orthonormal normal frame field of $M^{2}$: 
\begin{equation}
N_{1}=n(v);\text{ \ }N_{2}=-g^{\prime }(u)\,r(v)+f^{\prime }(u)\,e_{4}.
\label{C4}
\end{equation}%
Thus we obtain a positive orthonormal frame field $\{X,Y,N_{1},N_{2}\}$ of $%
M^{2}$. If we denote by $\kappa _{\alpha }(u)$ the curvature of the meridian
curve $\alpha (u)$, i.e. 
\begin{equation}
\kappa _{\alpha }(u)=f^{\prime }(u)\,g^{\prime \prime }(u)-g^{\prime
}(u)f^{\prime \prime }(u)={\frac{-f^{\prime \prime }(u)}{\sqrt{1-f^{\prime
2}(u)}}}.  \label{C5}
\end{equation}

Using (\ref{C31}) and (\ref{C4}) we can calculate the coefficients of the
second fundamental form of $X(u,v)$ as follows;%
\begin{eqnarray}
c_{11}^{1} &=&0,c_{22}^{1}=f(u)\kappa (v),  \notag \\
c_{12}^{1} &=&c_{12}^{2}=0,  \notag \\
c_{11}^{2} &=&\kappa _{\alpha }(u),  \label{C6} \\
c_{22}^{2} &=&f(u)g^{\prime }(u).  \notag
\end{eqnarray}

\begin{lemma}
Let $M^{2}$ be a meridian surface given with the surface patch (\ref{C2})
then 
\begin{equation}
A_{N_{_{1}}}=\left[ 
\begin{array}{ll}
0 & 0 \\ 
0 & {\frac{\kappa (v)}{f(u)}}%
\end{array}%
\right] ,\text{ }A_{N_{_{2}}}=\left[ 
\begin{array}{ll}
\kappa _{\alpha }(u) & 0 \\ 
0 & {\frac{g^{\prime }(u)}{f(u)}}%
\end{array}%
\right] .  \label{C8}
\end{equation}
\end{lemma}

Further by the use of (\ref{B6}) and (\ref{B7}) with (\ref{C6}), the Gauss
curvature is given by 
\begin{equation}
K={\frac{\kappa _{\alpha }(u)g^{\prime }(u)}{f(u)}.}  \label{C9}
\end{equation}

and the mean curvature vector field of $M^{2}$ becomes%
\begin{equation}
\overrightarrow{H}={\frac{\kappa (v)}{2f(u)}N_{_{1}}+\frac{\kappa _{\alpha
}(u)f(u)+g^{\prime }(u)}{2f(u)}N_{_{2}}.}  \label{C10}
\end{equation}

From the equation (\ref{C10}), we get the mean curvature of $M^{2}$

\begin{equation}
H=\frac{1}{2f(u)}\sqrt{\kappa (v)^{2}+\left( \kappa _{\alpha
}(u)f(u)+g^{\prime }(u)\right) ^{2}}.  \label{C11}
\end{equation}

\section{Proof of the Main Theorem}

Let $M^{2}$ be meridian surface given with the surface patch (\ref{C2}).
Then differentiating $K$ and $H$ with respect to $u$ and $v$ one can get%
\begin{eqnarray*}
K_{v} &=&0,\text{ }K_{u}=-\frac{\left( f(u)f^{^{\prime \prime \prime
}}(u)-f^{^{\prime }}(u)f^{^{\prime \prime }}(u)\right) }{f(u)^{2}}, \\
H_{v} &=&\frac{\kappa (v)\kappa ^{^{\prime }}(v)}{2f(u)\sqrt{\kappa
(v)^{2}+\left( \kappa _{\alpha }(u)f(u)+g^{\prime }(u)\right) ^{2}}}.
\end{eqnarray*}

Suppose that $M^{2}$ is a Weingarten surface then by the use of equation (%
\ref{A3}), we get,%
\begin{equation}
\frac{-\kappa (v)\kappa ^{^{\prime }}(v)\left( f(u)f^{^{\prime \prime \prime
}}(u)-f^{^{\prime }}(u)f^{^{\prime \prime }}(u)\right) }{2f(u)^{3}\sqrt{%
\kappa (v)^{2}+\left( \kappa _{\alpha }(u)f(u)+g^{\prime }(u)\right) ^{2}}}%
=0.  \label{C13}
\end{equation}

Thus we distinguish the following cases:

Case I: $\kappa (v)=0;$

\text{Case }II$:$ $\kappa ^{^{\prime }}(v)=0;$

Case III$:$ $f(u)f^{^{\prime \prime \prime }}(u)-f^{^{\prime
}}(u)f^{^{\prime \prime }}(u)=0.$

Let us consider \ these in turn;

Case I: Suppose $\kappa (v)=0$, i.e. the curve $c$ is a great circle on $%
S^{2}(1)$. In this case $N_{1}$ = const, and $M^{2}$ is a planar surface
lying in the constant $3$-dimensional space spanned by $\{X,Y,N_{2}\}$.
Particularly, if in addition $\kappa _{\alpha }(u)=0$, i.e. the meridian
curve lies on a straight line, then $M^{2}$ is a developable surface in the
3-dimensional space span $\{X,Y,N_{2}\}$ \cite{GM}.

Case II: Suppose $\kappa ^{^{\prime }}(v)=0.$ This implies that $\kappa (v)$
is nonzero constant. Then we have the following subcases;

Case II(a): $\kappa _{\alpha }(u)=0.$ In this case $c$ is a circle on $%
S^{2}(1)$, then $M^{2}$ is a developable ruled surface in a $3$-dimensional
Euclidean space $\mathbb{R}^{3}$.

Case II(b): $\kappa _{\alpha }(u)$ is nonzero constant$.$ In this case we
obtain the following ordinary differential equation.

\begin{equation}
{\frac{-f^{\prime \prime }(u)}{\sqrt{1-f^{\prime 2}(u)}}=a.}  \label{C13*}
\end{equation}%
Thus, the following expression is obtained from the solution of the
differential equation (\ref{C13*})%
\begin{equation*}
f(u)=\frac{\cos \left( au+ac_{1}\right) }{a}+c_{2}.
\end{equation*}%
Further, using the condition $(f^{\prime }(u))^{2}+(g^{\prime }(u))^{2}=1$
we get 
\begin{equation*}
g(u)=\frac{2\left( \sin \left( au+ac_{1}\right) -1\right) \sqrt{1+\sin
\left( au+ac_{1}\right) }}{\cos \left( au+ac_{1}\right) }.
\end{equation*}

Case III: $\ $Suppose $f(u)f^{^{\prime \prime \prime }}(u)-f^{^{\prime
}}(u)f^{^{\prime \prime }}(u)=0.$ Then we have the following subcases;

Case III(a): $f^{^{\prime \prime }}(u)=0.$ This implies that $\kappa
_{\alpha }(u)=K=0$, i.e. the meridian curve is part of a straight line and $%
M^{2}$ is a developable ruled surface. If in addition $\kappa (v)\neq $
const, i.e. $c$ is not a circle on $S^{2}(1)$, then $M^{2}$ is a developable
ruled surface in $\mathbb{E}^{4}$ (\cite{GM}).

Case III(b): $f^{^{\prime \prime }}(u)\neq 0.$ In this case we obtain the
following ordinary differential equation.%
\begin{equation}
f(u)f^{^{\prime \prime \prime }}(u)-f^{^{\prime }}(u)f^{^{\prime \prime
}}(u)=0  \label{C13**}
\end{equation}

Thus, the following expression is obtained from the solution of the
differential equation (\ref{C13**})

\begin{equation*}
f(u)=\pm \frac{1}{2}\sqrt{-\frac{a}{\left( e^{\frac{u}{b}}\right) ^{2}\left(
e^{\frac{c}{b}}\right) ^{2}}}\left( \left( e^{\frac{u}{b}}\right) ^{2}\left(
e^{\frac{c}{b}}\right) ^{2}+1\right) .
\end{equation*}

Further, using the condition $(f^{\prime }(u))^{2}+(g^{\prime }(u))^{2}=1$
one can get 
\begin{equation*}
g(u)=\pm \frac{1}{2}\sqrt{\frac{4b^{2}\left( e^{\frac{u}{b}}\right)
^{2}\left( e^{\frac{c}{b}}\right) ^{2}+a\left( e^{\frac{u}{b}}\right)
^{4}\left( e^{\frac{c}{b}}\right) ^{4}-2a\left( e^{\frac{u}{b}}\right)
^{2}\left( e^{\frac{c}{b}}\right) ^{2}+a}{b^{2}\left( e^{\frac{u}{b}}\right)
^{2}\left( e^{\frac{c}{b}}\right) ^{2}}}
\end{equation*}%
where $a,$ $b,$ $c,$ $c_{1},$ $c_{2}$ are real constants. This completes the
proof of the theorem.

\end{document}